\documentclass{amsart}
\usepackage{amscd}
\usepackage{graphicx}

\newcommand{\bbC}{{\mathbb C}}

\newcommand{\cH}{{\mathcal H}}
\newcommand{\bbP}{{\mathbb P}}
\newcommand{\bbR}{{\mathbb R}}
\newcommand{\bbZ}{{\mathbb Z}}
\newtheorem{lemma}{Lemma}[section]
\newtheorem{prop}[lemma]{Proposition}
\newtheorem{thm}{Theorem}
\newtheorem{other}[lemma]{Theorem}

\newtheorem{cor}[lemma]{Corollary}
\newtheorem{defn}{Definition}

\begin{document}
\title{Closed geodesics on incomplete surfaces.}

\thanks{The authors would like to acknowledge the support of the 
Australian Research Council}
\author{Paul Norbury}
\address{Department of Mathematics and Statistics\\
University of Melbourne\\Australia 3010}
\email{pnorbury@ms.unimelb.edu.au}
\author{J. Hyam Rubinstein}
\address{Department of Mathematics and Statistics\\
University of Melbourne\\Australia 3010}
\email{rubin@ms.unimelb.edu.au}

\begin{abstract}
    We consider the problem of finding embedded closed geodesics on
    the two-sphere with an incomplete metric defined outside a point. 
    Various techniques including curve shortening methods are used.
\end{abstract}

\maketitle

\section{Introduction.}
The existence of closed geodesics on a Riemannian surface often
depends only on the topology of the surface.  Having proven existence,
such a result usually falls short of giving basic properties like
whether the geodesics are embedded, or where the geodesics are located
on the surface.  A simple case of the latter is whether a geodesic
contains or avoids a given point on the surface.

The general problem addressed in this paper is the existence of
embedded closed geodesics on incomplete Riemannian surfaces.  We
require such a geodesic to {\em avoid the incomplete points.} Any
incomplete point is required to be contained in neighbourhoods with
arbitrarily small radius and area.  We loosely say the metric is
defined outside a finite set of points.

We will consider the situation when a surface (in fact a two-sphere)
contains exactly one incomplete point.  Topologically, one can think
of those embedded closed geodesics that miss the incomplete point as
the complement of the space of based loops.  That is 
\[\{{\rm based\ loops}\}\subset\{{\rm all\ loops}\}\] 
and we wish to ask if all of the critical points of the energy
functional are contained inside the smaller space.  In finite
dimensions one can have $S^2\subset S^3$ with all critical points of a
Morse function on $S^3$ contained inside $S^2$, whereas for
$\bbR\bbP^2\subset\bbR\bbP^3$ any Morse function on $\bbR\bbP^3$ must
have a critical point outside of $\bbR\bbP^2$ by a parity argument. 
On the space of loops, such an argument cannot go through without
further assumptions, since there is an example of an incomplete
two-sphere that does not possess an embedded closed geodesic avoiding
the incomplete point.  We give this example in
Section~\ref{sec:countex}.

There are two main approaches to proving the existence of closed
geodesics on Riemannian surfaces.  The first uses an infinite
dimensional version of Morse theory.  In finite dimensions, Morse
theory enables one to use the topology of a manifold to deduce
information about critical points of a smooth function on the
manifold.  By applying Morse theory to the length or energy functional
on the loop space of a Riemannian manifold, one aims to use the
topology of the loop space to deduce information about the existence
of closed geodesics.  The gradient flow, or path of steepest descent,
of the length function is not well-defined so one usually looks to
other curve-shortening methods to get existence and properties of
closed geodesics.  This is the main approach we take in this paper.

The study of periodic orbits of a Hamiltonian system on the tangent
bundle associated to the geodesic equation is the second main approach
for finding closed geodesics.  Various tools are used to study
Hamiltonian systems, such as methods of ordinary differential
equations (ODEs), Poincare sections and Floer homology.  Although we
will not use Floer homology here, it is worth pointing out that it
involves the study of an action functional on the loop space of the
tangent bundle using another infinite dimensional version of Morse
theory, where one replaces the gradient flow equation with a natural
holomorphic-like partial differential equation.  In
Section~\ref{sec:ode} we use ODE methods to study the Hamiltonian
system.

A model example of an incomplete two-sphere is an ellipsoid with
exactly three embedded closed geodesics, each given by the
intersection of the ellipsoid with a coordinate plane, with a point of
intersection of two of the geodesics removed.  Then there is precisely
one embedded closed geodesic on the incomplete surface.  More
generally, the three geodesics theorem gives the existence of three
embedded closed geodesics on a smooth Riemannian two-sphere.  It uses
the fact that the space of unparametrised loops on the two-sphere,
modulo identification of all constant loops, is homotopy equivalent to
$\bbR\bbP^3$.  Essentially, an embedded closed geodesic is associated
to each primitive homology class of $\bbR\bbP^3$, with the index of
the geodesic given by the dimension of the homology class.  Inside the
space of unparametrised loops is the set of loops containing a given
point on the two-sphere.  This set is homotopy equivalent to
$\bbR\bbP^2$ so one might expect to get only two embedded closed
geodesics containing the point, of index one and two.  (The index of a
closed geodesic might decrease when we hold fixed a point on the
surface.)  If one removes a point from a Riemannian two-sphere,
leaving an incomplete Riemannian metric on a surface, then one would
expect at least one geodesic to survive---avoid the incomplete point. 
Again, the counterexample of Section~\ref{sec:countex} means that one
needs to be careful with such an argument.

The Gauss-Bonnet theorem---that the integral of the Gaussian curvature
over a smooth closed Riemannian surface is $2\pi$ times the Euler
characteristic---does not apply to incomplete metrics.  One can use a
local version of the Gauss-Bonnet theorem applied to a neighbourhood
of a point, that includes the geodesic curvature of the boundary of
the neighbourhood, to measure the incompleteness at a point.  See
(\ref{eq:gb}) in Section~\ref{sec:circsym}.

The counterexample of Section~\ref{sec:countex} shows that it is
difficult in general to prove the existence of an embedded closed
geodesic on an incomplete surface.  However there is a circumstance 
that arises naturally where existence can be proven.  Consider an 
incomplete metric on a surface minus some points, that satisfies the 
two extra conditions:

(i) it has positive Gaussian curvature; and

(ii) the Gauss-Bonnet theorem fails in a strong way---the integral of
the Gaussian curvature over the surface is infinite.

An example of a metric satisfying such conditions, defined in a
neighbourhood of an incomplete point, is
\begin{equation}  \label{eq:modmet}
    ds^2=dr^2+rd\theta^2
\end{equation}
where $r$ is the distance from the point.  The curvature of this 
metric is $1/(2r^2)$.

\begin{thm}   \label{th:main}
    There exists an embedded closed geodesic on any two-sphere minus a
    point equipped with an incomplete positive curvature metric
    asymptotic to (\ref{eq:modmet}) and satisfying either of:
    
    (i) it admits a finite cyclic symmetry of order greater than 2;
    
    (ii) its shortest geodesic is long.
\end{thm}
The technical assumptions (i) or (ii) are necessary since the
counterexample in Section~\ref{sec:countex} consists of an incomplete
positive curvature metric asymptotic to (\ref{eq:modmet}).  Condition
(ii) is a condition on the ``length'' of a homology class of loops on
the two-sphere stated precisely in Definition~\ref{th:long}.  It can
be estimated in some cases.

The metrics in Theorem~\ref{th:main} arise naturally from
higher-dimensional problems.  Two applications are given in the next
two theorems.  In Section~\ref{sec:ode} we use ODE methods to prove
the existence of embedded closed geodesics in very particular
circumstances.  This allows us to show that condition (ii) of
Theorem~\ref{th:main} holds on the incomplete two-spheres
corresponding to the round three-sphere and the Fubini-Study metric on
$\bbC\bbP^2$.  Thus, a sufficiently close metric still satisfies
condition (ii) of Theorem~\ref{th:main} and has positive Gaussian
curvature.  One could pose the theorems in another way if discrete
symmetries are present.

\begin{thm}  \label{th:sphere}
    A three-sphere equipped with a circle invariant metric
    sufficiently close to the round metric possesses an embedded
    minimal torus invariant under the circle action.
\end{thm}
Theorem~\ref{th:sphere} should be compared with the result of White
\cite{WhiSpa} who proves the existence of a minimal torus for any
metric on the three-sphere sufficiently close to the round metric,
though without knowing if the torus is invariant under a circle
action.  In fact, we have chosen the round three-sphere rather
arbitrarily as an application of Theorem~\ref{th:main}.  There are
circle invariant metrics on the three-sphere with Ricci curvature not
everywhere positive that give rise to a positive curvature incomplete
metric on the two-sphere satisfying condition (ii) of
Theorem~\ref{th:main}, and thus possessing an embedded minimal torus.

The following theorem is true for a larger class of toric varieties
than $\bbC\bbP^2$, for example $S^2\times S^2$ and $\bbC\bbP^2$ blown
up twice.  In fact, it is believed that no assumptions of
``sufficiently close'' or discrete symmetries are necessary to get an
embedded minimal 3-torus.  We suspect that the counterexample in
Section~\ref{sec:countex} cannot be adjusted to apply to toric 
surfaces.
\begin{thm}  \label{th:cp2}
    A toric Kahler metric on $\bbC\bbP^2$ sufficiently close to the
    Fubini-Study metric possesses an embedded minimal 3-torus
    invariant under the torus action.
\end{thm}

Previously, closed geodesics on surfaces with incomplete metrics have
been studied using methods from ODEs \cite{HHTExi,HLaMin}.  They arose
in the context of minimal surfaces in higher dimensions, as in the
applications above.  The methods require the surfaces to be highly
symmetric.  In Section~\ref{sec:ode} we give such an argument for
highly symmetric Riemannian manifolds as a first step to proving a
more general existence result when the metric is slightly perturbed so
that it loses it symmetries.  Section~\ref{sec:flow} contains the
proof of Theorem~\ref{th:main}.  The applications are described in
detail in the final few sections.

We expect the existence of an embedded closed geodesic to be true in
greater generality, perhaps with the assumptions of
Theorem~\ref{th:main} minus the requirement that the Gaussian
curvature be positive away from the incomplete point.  In
Section~\ref{sec:minimax} we describe the minimax technique and the
intuition it brings to the problem.  A {\em sweepout} of a manifold is
a foliation (with singularities) of the manifold by a path of
codimension 1 submanifolds that degenerate to lower dimensional
submanifolds at both ends of the path.  Well-known examples are
sweepouts of a surface by circles, and sweepouts of the three-sphere
by tori that degenerate to the two core circles at each end.  Another
example is a sweepout of $\bbC\bbP^2$ by $T^3$s that degenerate to a
two-torus at one end and the union of three lines at the other, i.e.
take three lines $L_i\subset\bbC\bbP^2$, $i=1,2,3$ that do not
intersect at a common point.  The boundary of a disk neighbourhood of
the three lines is homeomorphic to $T^3$, and the complement of the
disk neighbourhood is homeomorphic to a neighbourhood of a
homologically trivial $T^2\subset\bbC\bbP^2$.

A sweepout has a maximal volume leaf.  The {\em minimax} of a family
of sweepouts is the infimum of the volumes of the maximal volume
leaves.  If the family is large enough the minimax is minimal, and
very often the minimax (or perhaps a piece of it) is a minimal
submanifold.

On a two-sphere minus a point with an incomplete Riemannian metric
there is a natural family of sweepouts by loops.  The loops degenerate
to the incomplete point at one end of the sweepout and to an interior
point at the other end.  Although the maximum length loop in each
sweepout lies inside the two-sphere minus a point, unfortunately the
minimax of a family of such sweepouts may not avoid the incomplete
point.  The maximum length loops may gradually move towards the
incomplete point.  In Section~\ref{sec:change} we give two examples of
sequences of sweepouts over a two-sphere minus a point with minimax
not contained in the two-sphere minus a point.  They correspond to a
sequence of sweepouts of the round three-sphere by tori and
$\bbC\bbP^2$ by three-tori, with a change of topology in the
minimax---rather than obtaining the Clifford torus or the minimal
$T^3\subset\bbC\bbP^2$ the minimax is a minimal two-sphere,
respectively minimal three-sphere.

Although we insist that geodesics on an incomplete surface avoid 
incomplete points, one can drop this condition since shortest paths 
still exist on the compact manifold with a singular metric.  In 
\cite{GHaClo} it is proven that when the incomplete surface is an 
orbifold there exists a closed geodesic.  In the case of a two-sphere 
with one orbifold point, the closed geodesic contains the orbifold 
point.

In the following sections we prove the existence of embedded closed
geodesics on the two-sphere minus a point equipped with an incomplete
metric in increasingly more general circumstances.  We begin with the
circle symmetric case where the corresponding Hamiltonian system is
integrable and closed geodesics are essentially understood.  We move
on to a class of metrics with discrete symmetries.  Next we add the
assumption of positive Gaussian curvature with infinite integral, and
consider this case with and without finite symmetries.  Finally we 
study minimax of sweepouts and larger families of loops.

\section{ODE methods}  \label{sec:ode}
It is often convenient and natural to represent an incomplete point on
a surface as a circle boundary along which the metric is degenerate. 
Polar coordinates $(r,\theta)$ give such a representation around an 
incomplete point as in (\ref{eq:modmet}) or around a complete point.  
We shall take this viewpoint to represent the two-sphere as a disk 
with incomplete point given by its boundary.  

\subsection{Circle symmetric metrics}  \label{sec:circsym}
The simplest family of incomplete metrics on the two-sphere is
\begin{equation}  \label{eq:simple}
    ds^2=dr^2+f(r)^2d\theta^2
\end{equation}
where $r\in[0,1]$ and $\theta\in[0,2\pi]$, and $f$ vanishes at $r=0$
and $r=1$.  The case $f(r)=r\sqrt{1-r^2}$ was studied in \cite{HLaMin}
(using different coordinates) giving rise to minimal tori in the round
three-sphere.  

We choose $f$ to look like $r$ near $0$, precisely $f'(0)=1$, so that
the metric is complete there.  This follows from the local
Gauss-Bonnet formula
\begin{equation} \label{eq:gb}
    \int_{\Omega}Kda+\int_{\partial\Omega}kds=2\pi
\end{equation}
applied to the disk $r\leq\epsilon$, where $K=-f''/f$ is the 
Gaussian curvature, $k=f'/f$ is the geodesic curvature of the 
curve $r\equiv\epsilon$, $da=fdrd\theta$ and $ds=fd\theta$.  After 
integrating out the $\theta$ terms (\ref{eq:gb}) becomes
\[ -\int_0^{\epsilon}f''dr+f'(\epsilon)=1\]
and the left hand side is $f'(0)$.  We allow $f$ to vanish in any way
at $r=1$, although we say the metric is incomplete at the `point'
$r=1$ whether or not it can be completed.  For example if it looks
like $1-r$, respectively $\sqrt{1-r}$, then in the former case it can
be completed, and in the latter it cannot.

Geodesics for (\ref{eq:simple}) parametrised by arc-length satisfy
\begin{equation}  \label{eq:geodsimp}
    f(r)^2
    \dot{\theta}=c,\ \ \ \ \ \dot{r}^2+f(r)^2\dot{\theta}^2=1
\end{equation}    
for some constant $c$.  It immediately follows that $f(r)\geq c^2$ so
geodesics either meet the incomplete point when $c=0$, or they are
confined to remain a bounded distance from the incomplete point since
$c\neq 0$.  In the latter case, by the symmetry of the metric each
geodesic is periodic in $\theta$, measured between two closest points
to $r=0$ or $r=1$, and when the period is a rational multiple of
$2\pi$, the geodesic is closed.  

Each critical point of $f$ in the interior of $(0,1)$ corresponds to
an embedded closed geodesic, so $f_r(r_0)=0$ implies $r\equiv r_0$ is
a geodesic.  This occurs at least once, at the maximum of $f(r)$, and
precisely once when the Gaussian curvature is positive, since then
$f(r)$ is convex and its maximum is the unique stationary point in
$(0,1)$.  

In the example $f(r)=r\sqrt{1-r^2}$ where embedded closed geodesics
correspond to embedded minimal tori in the three-sphere, the
incomplete metric has positive Gaussian curvature with $\int
K=\infty$.  The only embedded closed geodesic comes from the maximum
of $f(r)$ and it corresponds to the Clifford torus.  Nevertheless, in
general when a metric can not be completed, i.e. $\int K\neq 4\pi$,
there may be another embedded closed geodesic as shown in the
following proposition.

Let $\gamma$ be the embedded closed geodesic corresponding to the
maximum of $f(r)$. 
\begin{prop}
    The metric (\ref{eq:simple}) has at least two embedded closed
    geodesics if
    \[\int K<4\pi\ \ \ {\rm and}\ \ {\rm index\ }\gamma>1\]
    \[{\rm or}\ \ \ \ \ \int K>4\pi\ \ \ {\rm and}\ \ \ {\rm index\
    }\gamma=1.\]
\end{prop}
\begin{proof}
    Specify each geodesic by the minimum value $c$ that $f$ takes on
    the geodesic.  As described above, each geodesic is periodic in
    $\theta$, and the period $\Omega_c$ is a rational multiple of
    $2\pi$ for closed geodesics.  {\em Embedded} closed geodesics have
    period $2\pi/n$ for a positive integer $n$.  If 
    \begin{equation}  \label{eq:perineq}
	\lim_{c\rightarrow 0}\Omega_c>2\pi>
	\lim_{c\rightarrow c_{\rm crit}}\Omega_c
    \end{equation}
    then by continuity of the period and the intermediate value
    theorem there must be a geodesic with $r$ non-constant and period
    $2\pi$.  This gives an embedded closed geodesic other than
    $\gamma$.  
    
    For the moment assume the following two limits which we prove
    below:
    \begin{equation}   \label{eq:limper}
	\lim_{c\rightarrow 0}\Omega_c=\pi-\pi/f'(1),\ \ \ 
	\lim_{c\rightarrow c_{\rm crit}}\Omega_c
	=2\pi/\sqrt{-f''(r_{\rm crit})c_{\rm crit}}
    \end{equation}    
    where $f'(1)\in[-\infty,0)$ and $f(r_{\rm crit})=c_{\rm crit}$ is the 
    maximum of $f$.  Now, $L=2\pi c_{\rm crit}$, 
    $K|_L=-f''(r_{\rm crit})/c_{\rm crit}$ and
    \[\int Kda=\int_0^{2\pi}\int_0^1\frac{-f''(r)}{f}f 
    drd\theta=-2\pi f'(r)|^1_0=2\pi-2\pi f'(1)\]
    so
    \[\lim_{c\rightarrow 0}\Omega_c>2\pi\ \  (<2\pi)\ 
    \Leftrightarrow\ \int K<4\pi\ \ (>4\pi).\]
    
    The Gaussian curvature is constant along $\gamma$ so the Rauch 
    comparison theorem tells us the exact distance between conjugate 
    points along $\gamma$---it is $\pi/\sqrt{K}$ where we have abused 
    notation and put $K=-f''(r_{\rm crit})/c_{\rm crit}$ for the value 
    of the Gaussian curvature on $\gamma$.  Then
    \[ L\sqrt{K}=2\pi c_{\rm crit}\times 
    \sqrt{-f''(r_{\rm crit})/c_{\rm crit}}=
    2\pi\sqrt{-f''(r_{\rm crit})/c_{\rm crit}}.\]
    Thus
    \[\lim_{c\rightarrow c_{\rm crit}}\Omega_c<2\pi\ \  (>2\pi)\
    \Leftrightarrow\ L\sqrt{K}|>2\pi\ \  (<2\pi).\]
    Finally
    \begin{equation}   \label{eq:index}
	{\rm index\ }\gamma>1\Leftrightarrow L\sqrt{K}>2\pi
    \end{equation}
    since if $L\sqrt{K}>2\pi$ then $\gamma$ can be cut into two arcs
    of length greater than $\pi/\sqrt{K}$ so there are pairs of
    conjugate points on the arcs and the index is at least two. 
    Conversely, if index~$\gamma>1$ then $\gamma$ must contain at
    least two disjoint arcs of length greater than $\pi/\sqrt{K}$. 
    Now, the index of the geodesic is at least 1 since the geodesic
    can be decreased in size (take nearby constant $r$ loops) so the
    the hypotheses of the proposition imply (\ref{eq:perineq}) with
    its inequalities as written, or reversed.\\
    
    It remains to prove (\ref{eq:limper}).  From (\ref{eq:geodsimp}),
    \[\frac{d\theta}{dr}=\frac{\dot{\theta}}{\dot{r}} 
    =\frac{c}{f(r)\sqrt{f(r)^2-c^2}}\]
    so the period of the geodesic is
    \[\Omega_c=2\int_{r_1}^{r_2}\frac{c}{f(r)\sqrt{f(r)^2-c^2}}dr\]
    where $f(r_1)=c$ and $f(r_2)=c$ is the next time along the 
    geodesic $f$ is $c$.
    
    The limit $\lim_{c\rightarrow 0}\Omega_c$ depends only on the 
    local behaviour of $f$ near $0$ and $1$ since for any $\epsilon>0$
    \[ \lim_{c\rightarrow 0}\int_{\epsilon}^{1-\epsilon}
    \frac{1}{f(r)\sqrt{f(r)^2-c^2}}dr=\int_{\epsilon}^{1-\epsilon}
    \frac{dr}{f(r)^2}=<\infty\]
    and hence it is annihilated by multiplication by $c$ when 
    $c\rightarrow 0$.  Thus
    \[\lim_{c\rightarrow 0}\Omega_c=2\lim_{c\rightarrow 0}\left(
    \int_{r_1}^{\epsilon}\frac{c}{f(r)\sqrt{f(r)^2-c^2}}dr
    +\int_{1-\epsilon}^{r_2}\frac{c}{f(r)\sqrt{f(r)^2-c^2}}dr\right).\]
    
    Choose $\alpha>f'(0)$ and $\epsilon$ small enough so that for each
    $r_1\in(0,\epsilon]$, $f(r)$ lies below the line
    $y=c+\alpha(r-r_1)$ for $r\in[r_1,\epsilon]$.  This gives a lower
    bound for the limit
    \begin{eqnarray*}
	2\lim_{c\rightarrow 0}\int_{r_1}^{\epsilon}
	\frac{c}{f(r)\sqrt{f(r)^2-c^2}}dr&\geq&
	2\lim_{c\rightarrow 0}\int_c^{c_1}
	\frac{c}{\alpha y\sqrt{y^2-c^2}}dy\\
	&=&2\lim_{c\rightarrow 0}\frac{1}{\alpha}\arctan 
	{\sqrt{(c_1/c)^2-1}}\\
	&=&\pi/\alpha
    \end{eqnarray*}
    where $c_1=c+\alpha(\epsilon-r_1)$.  Similarly, we can choose
    $\alpha<f'(0)$ so that $f(r)$ lies above a family of lines with
    slope $\alpha$ to get an upper bound for the limit.  The same
    argument works near $r=1$, using $f'(1)$.  If $f'(1)=-\infty$ then
    we simply get an upper bound, with the lower bound of 0 automatic. 
    Thus
    \[\lim_{c\rightarrow 0}\Omega_c=\pi/f'(0)+\pi/|f'(1)|=
    \pi-\pi/f'(1).\]
    
    The limit near the critical point, which we may translate to
    $r=0$, is a local quantity.  It can be simplified as follows.
    \[\lim_{c\rightarrow c_{\rm crit}}\Omega_c=
    \lim_{c\rightarrow c_{\rm crit}}2\int_{r_1}^{r_2}\frac{c} 
    {f(r)\sqrt{f(r)^2-c^2}}dr=
    \lim_{c\rightarrow c_{\rm crit}}\frac{2}{\sqrt{2c_{\rm crit}}}
    \int_{r_1}^{r_2}\frac{dr}{\sqrt{f-c}}.\]
    If we change the parametrisation to $r=r(t)$ then the limit transforms as
    \[\lim_{c\rightarrow c_{\rm crit}}
    \int_{r_1}^{r_2}\frac{dr}{\sqrt{f-c}}=
    \lim_{c\rightarrow c_{\rm crit}}\int_{t_1}^{t_2}\frac{r'(t)dt}
    {\sqrt{f-c}}=
    \lim_{c\rightarrow c_{\rm crit}}\int_{t_1}^{t_2}\frac{r'(0)dt}
    {\sqrt{f-c}}\]
    and
    \[\frac{df^2}{dt^2}(0)=\frac{d^2f}{dr^2}(0)r'(0)^2+\frac{df}{dr}(0)r''(0)
    =\frac{df^2}{dr^2}(0)r'(0)^2.\]
    Thus if we choose $r'(0)=1$ then in the limit $f(r)$ is simply
    replaced by $f(r(t))$ and $f''(0)$ is well-defined.  By the Morse
    lemma, since the critical point is non-degenerate we can choose
    $r(t)$ so that $f=c_{\rm crit}-\lambda^2t^2$ where
    $\lambda^2=-f''(0)/2$.  The limit becomes
    \[\lim_{c\rightarrow c_{\rm crit}}\int_{t_1}^{t_2}\frac{dt}
    {\sqrt{c_{\rm crit}-c-\lambda^2t^2}}=\pi/\lambda.\]
    so \[\lim_{c\rightarrow c_{\rm crit}}
    \Omega_c=2\pi/\sqrt{-f''(r_{\rm crit})c_{\rm crit}}.\] 
\end{proof}
When the Gaussian curvature is positive and $\int K=\infty$, the case
we study more generally in Section~\ref{sec:flow}, one can replace the
hypothesis index~$\gamma=1$ by the stronger hypothesis $K<1$ on
$\gamma$, i.e. the function $f(r)$ is rather flat at its maximum. 
This follows by the convexity of $f(r)$ and $f'(0)=1$ which imply
$L<2\pi$, and together with the assumption $K<1$, gives
$L\sqrt{K}<2\pi$ so by (\ref{eq:index}) this implies index~$\gamma=1$.

\subsection{Polygons}  \label{sec:poly}
A more complicated class of incomplete metrics on the two-sphere arise
from toric geometry.  See Section~\ref{sec:toric}.  Using ODE methods,
we will analyse a specific example which possesses discrete
symmetries.

On the square $|x|<1$, $|y|<1$ define the metric
\begin{equation}  \label{eq:metprod}
  ds^2=(1-y^2)dx^2+(1-x^2)dx^2.
\end{equation}
This is an incomplete metric on the two-sphere, with the boundary the 
incomplete point.

The equations for the geodesic flow are 
\begin{eqnarray*}
    \ddot{x}&=&\frac{2y}{1-y^2}\dot{x}\dot{y}-\frac{x}{1-y^2}\dot{y}^2\\
    \ddot{y}&=&\frac{-y}{1-x^2}\dot{x}^2+\frac{2x}{1-x^2}\dot{x}\dot{y}
\end{eqnarray*}
or implicitly
\begin{equation}  \label{eq:geodimp}
    \frac{d^2y}{dx^2}=\frac{x}{1-y^2}\left(\frac{dy}{dx}\right)^3
    -\frac{2y}{1-y^2}\left(\frac{dy}{dx}\right)^2+\frac{2x}{1-x^2}
    \frac{dy}{dx}-\frac{y}{1-x^2}.
\end{equation}
A compact way to respresent the geodesic equations is through the
equivalent Hamiltonian system
\begin{equation}
    H=\frac{1}{2}\left(\frac{p_1^2}{1-q_2^2} 
    +\frac{p_2^2}{1-q_1^2}\right).
\end{equation}
We will not use the Hamiltonian formulation, except to point out that 
it appears to be non-integrable so we expect the geodesic flow to be 
complicated.

The remainder of this section is devoted to the proof that there is 
an embedded closed geodesic on the two-sphere with incomplete metric 
given by (\ref{eq:metprod}).

Consider the family of geodesics $\cH$ that meet the $y$-axis
horizontally---$\dot{y}=0$ at $x=0$.  The next few lemmas prove that
each geodesic in $\cH$ meets the line $x=y$ and at least one of these
meets $x=y$ orthogonally.  By symmetry, this geodesic extends to a
geodesic that is an embedded loop.  Furthermore, we get estimates on
the length and the index of the closed geodesic.

\begin{lemma}
    A geodesic that meets the $y$-axis twice while remaining inside
    the region $|y|\geq|x|$ must be the $y$-axis.
\end{lemma}
\begin{proof}
    A geodesic that meets the $y$-axis twice while $|y|\geq|x|$ must
    be tangent somewhere to a line $y=cx$, and it must lie to the
    $y$-axis side of the line (unless the geodesic is the $y$-axis.) 
    From (\ref{eq:geodimp}), when $y=cx$ and $dy/dx=c$ we have
    \begin{eqnarray*}
	\frac{d^2y}{dx^2}&=&\frac{x}{1-c^2x^2}c^3
    -\frac{2cx}{1-c^2x^2}c^2+\frac{2x}{1-x^2}c-\frac{cx}{1-x^2}\\
    &=&\frac{cx}{(1-x^2)(1-c^2x^2)}\left(-(1-x^2)c^2+(1-c^2x^2)\right)\\
    &=&\frac{y}{(1-x^2)(1-y^2)}(1-c^2)
    \end{eqnarray*}
    Since $|y|>|x|$, $|c|>1$ so $1-c^2<0$ and $d^2y/dx^2$ has the
    opposite sign of $y$.  But this means the geodesic and the
    $y$-axis are on opposite sides of the tangent line, which is a
    contradiction.
\end{proof}

\begin{lemma}
    Each geodesic in $\cH$ meets the lines $x=\pm y$ with $y$ monotone
    in $t$.
\end{lemma}
\begin{proof}
    By symmetry, without loss of generality we may assume $y>x>0$.  A
    geodesic that satisfies $\dot{y}\leq 0$ and $y>0$ somewhere, must
    continue to have $\dot{y}<0$ until $y<0$---if $\dot{y}$ gets too
    close to zero while $y>0$ then $\ddot{y}<0$ so $\dot{y}$ is sent
    back away from zero.  This is an easy consequence of one of the
    equations for the geodesic,
    $\ddot{y}=(-y\dot{x}^2+2x\dot{x}\dot{y})/(1-x^2)$.  If
    $y>\epsilon>0$ and $\dot{y}$ is small enough then $\ddot{y}<0$
    since $\dot{x}$ is bounded.
    
    Thus, a geodesic in $\cH$ must travel down and pass through $y=0$ 
    while $\dot{y}<0$.  From the previous lemma, it cannot pass 
    through the $y$-axis while $y>0$ so it must meet the line $y=x$.
\end{proof}

\begin{lemma}
    Any geodesic in $\cH$ that meets the lines $x=\pm y$ at an angle 
    less than $\pi/2$ is convex.
\end{lemma}
\begin{proof}
    Without loss of generality assume $y\geq x\geq 0$.  Suppose a
    geodesic in $\cH$ meets the line $x=y$ at an angle less than
    $\pi/2$, so its angle to the horizontal is greater than $\pi/2$. 
    In the first quadrant $dy/dx\leq 0$ implies that $d^2y/dx^2<0$
    since each term of (\ref{eq:geodimp}) is negative.  In particular
    a geodesic that meets the positive $y$-axis horizontally has
    $dy/dx$ negative and decreasing for small positive $x$.  Either
    this behaviour continues until $y<x$, proving the lemma, or
    $dy/dx\rightarrow -\infty$ and the geodesic becomes vertical.  The
    latter case does not occur, because after it is vertical
    $\dot{x}<0$ since $d^2x/dy^2<0$ by exchanging $x$ and $y$ in the
    argument earlier in this paragraph.  Furthermore, it cannot become
    vertical again while $x>0$ since $d^2x/dy^2<0$.  In particular, it
    could not have met the line $x=y$ at an angle less than $\pi/2$.
\end{proof}

\begin{prop}
    There is an embedded closed geodesic in the square.
\end{prop}
\begin{proof}
    We will prove that there is a geodesic that meets the lines $x=0$,
    $y=0$ and $y=\pm x$ orthogonally.
    
    Any geodesic in $\cH$ close enough to the $x$-axis is almost
    horizontal so it meets the lines $y=\pm x$ at an angle
    approximately equal to $\pi/4$.  Below we will show that any
    geodesic in $\cH$ close enough to the boundary meets the lines
    $y=\pm x$ at an angle greater than $\pi/2$.  By continuity of the
    angle and the intermediate value theorem, there is a geodesic in
    between that meets the lines $y=\pm x$ at an angle of $\pi/2$. 
    The extension of this geodesic past the line $y=x$ is obtained by
    reflecting it through $y=x$.  Reflect twice more to obtain an
    embedded closed geodesic.
    
    Geodesics close to the boundary curve sharply away from the
    boundary.  This is because the metric near the boundary (away from
    the vertices of the square) behaves like the model metric
    (\ref{eq:modmet}) where the geodesics are explicit.  Now assume
    $y>x>0$.  By the proof of the previous lemma, once the slope of
    the geodesic is steeper than $-1$, it must meet the line $y=x$ at 
    an angle greater than $\pi/2$ since it is either convex all the 
    way to the line $y=x$, and hence steeper there, or it becomes 
    vertical somewhere and hence must meet $y=x$ at an angle greater 
    than $3\pi/4$.
\end{proof}

\begin{figure}[ht]
	\scalebox{1}{\includegraphics[height=6cm]{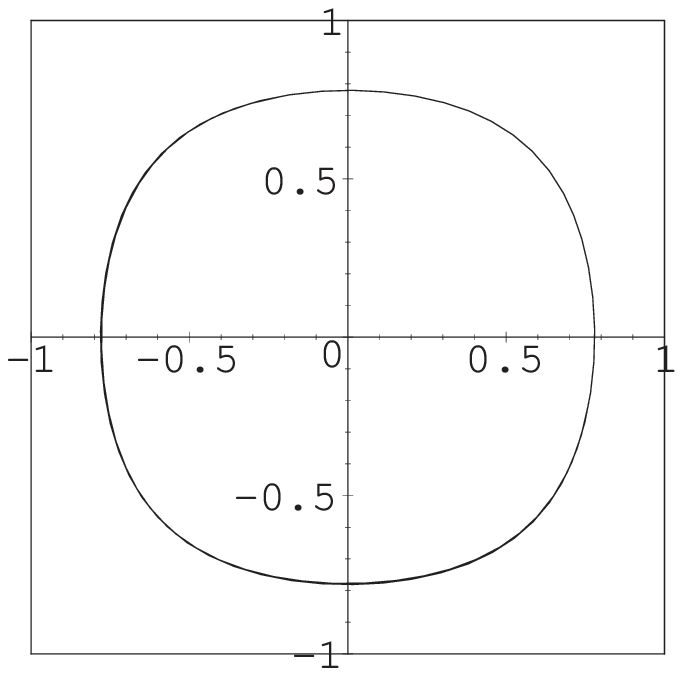}}
        \scalebox{1}{\includegraphics[height=6cm]{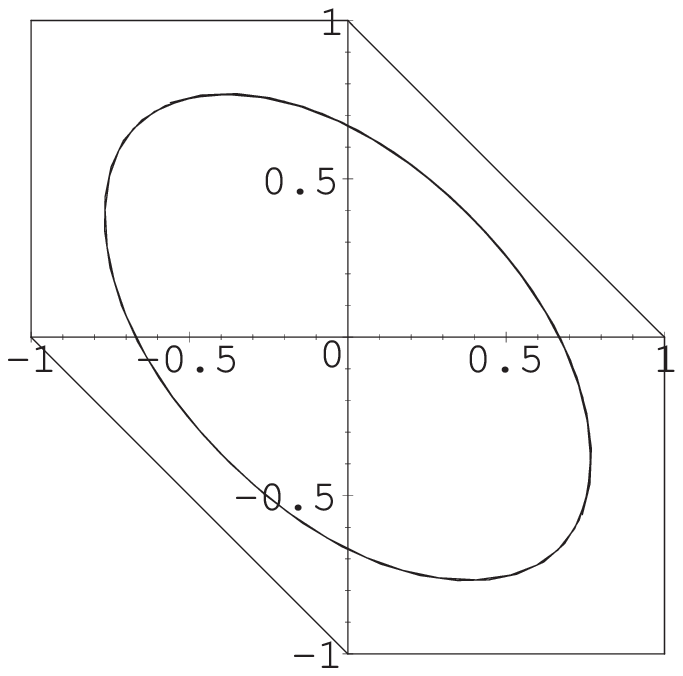}}
    \caption{Embedded closed geodesic.}
    \label{fig:minimal}
\end{figure}

The same argument applies to other regular polygons to give embedded
closed geodesics.  For example, the triangle and hexagon defined by
$abc=0$ for $(a,b,c)=(1+x,1+y,1-x-y)$, respectively
$(a,b,c)=(1-x^2,1-y^2,1-(x+y)^2$ equipped with the metric
$ds^2=(bcdx^2+acdy^2+abd(x+y)^2)/(a+b+c)$ are regular since their
affine symmetries, which preserve the metric, make up the whole
polygon symmetry group.  In Section~\ref{sec:toric} we will see that
the square, triangle and hexagon correspond to $S^2\times S^2$,
$\bbC\bbP^2$ and $\bbC\bbP^2$ blown up at three points, and the
embedded closed geodesic in each corresponds to an embedded minimal
$T^3$.  Figure~\ref{fig:minimal} shows the closed geodesics obtained
numerically using MAPLE. The geodesics lift to a minimal hypertorus in
$S^2\times S^2$, respectively $\bbC\bbP^2$ blown up at three points.

\section{Geodesic curvature flow.}  \label{sec:flow}

Gage \cite{GagDef} considered the flow of a loop $\gamma$ on a surface
of curvature $K>0$ 
\begin{equation}  \label{eq:flow}
    \frac{\partial\gamma}{\partial t}=\frac{k}{K}\nu
\end{equation}
in its normal direction $\nu$, where $k$ is the geodesic curvature of
$\gamma$.  The flow preserves embeddedness and the integral of $K$
over the area enclosed by the loop behaves well under the flow.  If
the integral is $2\pi$ then it is preserved by the flow.  Chou and Zhu
\cite{CZhCur} completed the proof of long time existence to prove that
there exists an embedded closed geodesic on a 2-sphere with positive
Gaussian curvature.

\subsection{Solutions in the neighbourhood of a point}
Around a complete point of the metric, any small loop flows to a point
in finite time whereas around an incomplete point of the metric, it
can take infinite time.  This can be seen explicitly for
$ds^2=dr^2+f(r)^2d\theta^2$ when $f$ is convex.  The curve
$r=$constant has
\[ k=-\frac{f'}{f},\ K=-\frac{f''}{f},\ \nu=(1,0)\]
so the flow is given by $dr/dt=f'/f''$ with solution
\[ f'(r)=f'(r_0)e^t.\]
In a neighbourhood of $r=0$, $0<f'(r_0)<f'(r)$ since $f$ is convex so 
when $f'(0)$ is finite, for example $f'(0)=1$ in the complete case, 
the loop converges to a point in finite time.

For the remainder of this section we will assume that $f'(0)$ is
infinite.  In this case the curve $r=$constant takes infinite time to
converge to a point.  For example, when the metric is the model metric
(\ref{eq:modmet}), so $ds^2=dr^2+rd\theta^2$ the flow is given by
\[ dr/dt=-2r\ \Rightarrow\ r=r_0e^{-2t}\] which shrinks exponentially.

When 
\begin{equation}  \label{eq:cornmet}
    ds^2=ydx^2+xdy^2
\end{equation}    
we do not find an exact solution of the flow however the levels sets
of $xy$ behave well under the flow.  The curve $xy=\epsilon$ has \[
k=\frac{x^2+y^2}{2xy(x+y)^{3/2}},\ K=\frac{x+y}{4x^2y^2},\
\nu=\frac{-1}{\sqrt{x+y}}(1,1)\] so the infinitesimal change in
$\epsilon=xy$ during the flow is \[
\frac{d\epsilon}{dt}=\nabla\epsilon\cdot(\dot{x},\dot{y})=
-2xy\frac{x^2+y^2}{(x+y)^2}.\]

Denote the geodesic flow of $xy=\epsilon_0$ by $(x(t),y(t))$.
We can use the fact that 
\begin{equation}  \label{eq:low}
    1\geq(x^2+y^2)/(x+y)^2\geq 1/2
\end{equation}
to prove that 
\[ \epsilon_0e^{-2t}\leq x(t)y(t)\leq\epsilon_0e^{-t}.\]
The flow $\gamma(t)=\{ xy=\epsilon_0e^{-2t}\}$ coincides with the
curve $xy=\epsilon_0$ at $t=0$ and by the first inequality in
(\ref{eq:low}), $\gamma(t)$ initially moves off more rapidly than the
geodesic flow $(x(t),y(t))$.  If $(x(t),y(t))$ was ever to catch
$\gamma(t)$, then it would have to meet it at a tangent such that its
geodesic curvature is bounded below by the geodesic curvature of
$\gamma(t)$.  But then $\gamma(t)$ is moving faster which contradicts
the fact that it was ever caught by $(x(t),y(t))$.  The same argument
shows that the curve $\{ xy=\epsilon_0e^{-t}\}$ never catches
$(x(t),y(t))$ when $t>0$.  Thus the curve takes infinite time to reach
the boundary and it moves exponentially fast there.

Any loop in a neighbourhood of the incomplete point (that does or does
not contain the incomplete point) flows toward the incomplete point in
infinite time.  This is because the level sets $r=$constant (or
$xy=$constant in the second example) form barriers for the loop,
forcing it to remain between two level sets that move exponentially fast
toward the incomplete point.

A loop with $\int K<2\pi$ does not actually reach the incomplete
point, instead it shrinks to a point (away from the incomplete point)
in finite time.  This is because we must have $\int_{\gamma}k>0$.  The
flow satisfies \[ \frac{d\int_{\gamma}k}{dt}=\int_{\gamma}k\] so
$\int_{\gamma}k$ increases to $2\pi$ in finite time.  At that time the
flow stops since $\int K=0$ inside the loop.  This occurs away from
the incomplete point due to an exponentially shrinking barrier so we
can apply the results of \cite{CZhCur} that if the flow exists only
for finite time then the loop flows to a point.

Before we can understand the behaviour of loops with $\int K=2\pi$ we
must first consider scale invariant solutions.  Consider solutions to
(\ref{eq:flow}) that are scale invariant under the map
\begin{equation}   \label{eq:rescale}
    (r,\theta)\mapsto(\lambda(t)^2r,\lambda(t)\theta) 
\end{equation}    
where $\lambda(t)=e^{ct}$ determines the speed of a solution.  (If we
choose $c=0$ then this produces a geodesic.)  In fact, we choose
$c=1$, in order to scale at the same rate as the flow of the curves
$r=$constant, or in other words, so that $r=$constant is a scale
invariant solution.  Lift to the cover with $\theta\in\bbR$.  Scale
invariant solutions are static solutions of the equation
\begin{equation}  \label{eq:scaleflow}
    \frac{dr}{dt}=4r^2\ddot{r}-r\theta\dot{r}\dot{\theta},\ \ \ 
    \frac{d\theta}{dt}=4r^2\ddot{\theta}+2r\dot{r}\dot{\theta}
    +\theta\dot{r}^2,\ \ \ 
    \dot{r}^2+r\dot{\theta}^2=1.
\end{equation}
The variable $t$ is the flow parameter and the dotted derivatives are
taken with respect to the arc length $s$.  The ODE obtained by setting
$dr/dt=0=d\theta/dt$ in (\ref{eq:scaleflow}) has a unique solution for
each set of initial data.  

\begin{figure}[ht]
    \scalebox{1}{\includegraphics[height=6cm]{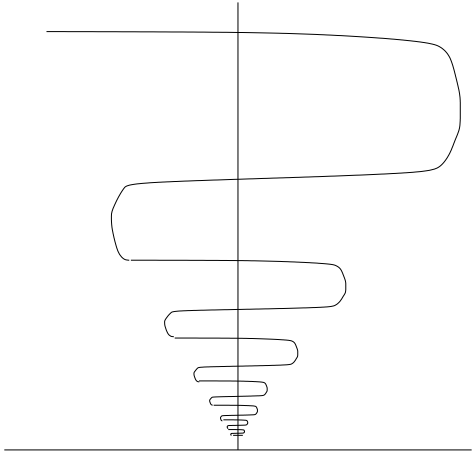}}
    \caption{Scale invariant solution $\xi_{a,b}$ acts as a barrier.}
    \label{fig:scale}
\end{figure}
Consider the ODE given by the static solution of \ref{eq:scaleflow}. 
Associate to any pair $(a,b)\in\bbR^+\times\bbR$ the unique solution
$\xi_{a,b}$ of this ODE with initial conditions at $s=0$ \[ r=a,\
\theta=b,\ \dot{r}=1,\ \dot{\theta}=0.\]
\begin{prop}
    The static solution $\xi_{a,b}$ reaches $\theta=0$ for $s>0$ and
    $s<0$ at angles $\alpha$ and $\beta$, and over the interior of the
    region bounded by the solution curve and the line $\theta=0$,
    $\int K=\alpha+\beta$.
\end{prop}
\begin{proof}
    The solution $\xi_{a,b}$ reaches $\theta=0$ for forward and
    backward $s$ since $\ddot{r}=\dot{r}\dot{\theta}\theta/(4r)$
    ensures that while $\theta>0$ (say $b>0$), $\dot{r}$ decreases as
    $s$ increases above 0 or decreases below 0.  Furthermore,
    $\dot{r}$ remains positive since if $\dot{r}=0$ then the solution
    must coincide with the solution $r=$constant which has no vertical
    part.  By $\dot{r}^2+r\dot{\theta}^2=1$, $|\dot{\theta}|$
    increases so the curve reaches $\theta=0$.
    
    The second claim $\int K=\alpha+\beta$ is equivalent to the
    statement that $\int k=0$ along $\xi_{a,b}$ running from
    $\theta=0$ to $\theta=0$ since the Gauss-Bonnet formula gives
    $\int K+\int k+\pi-\alpha+\pi-\beta=2\pi$.  Now
\begin{eqnarray*}
    k&=&r^{-1/2}(r\ddot{r}\dot{\theta}=r\dot{\theta}^3/2
    -r\ddot{\theta}\dot{r}-\dot{r}^2 \dot{\theta})\\
    &=&r^{-1/2}(\theta\dot{r}\dot{\theta}^2/4-r\dot{\theta}^3/2
    -\dot{r}^2\dot{\theta}/2+\theta\dot{r}^3/(4r))\\
    &=&r^{-1/2}(\dot{r}\theta/(4r)-\dot{\theta}/2)\\
    &=&\dot{((-1/2)r^{-1/2}\theta)}
\end{eqnarray*}
    where the first two equations of (\ref{eq:scaleflow}) were used to
    go from line 1 to line 2 and the third equation of
    (\ref{eq:scaleflow}) was used to go from line 2 to line 3.  Thus
    \[\int_{\gamma}kds=(-1/2)r^{-1/2}(\theta_2-\theta_1)\] but each
    $\theta_i=0$ so $\int_{\gamma}kds=0$ as claimed.  
\end{proof}
    
The shape of a solution is given in Figure~\ref{fig:scale}.  If
$b>>0$, then the solution becomes almost horizontal as $\theta$
approaches 0, i.e. the two angles $\alpha$ and $\beta$ of the solution
curve with $\theta=0$ are approximately $\pi/2$.  This is because
$\dot{r}$ must remain positive, $\alpha$ and $\beta$ are both less
than $\pi/2$ but then $\int K<\pi$ over the region bounded by the
solution curve and $\theta=0$.  In order to maintain such a small
integral when $b>>0$, the solution curve may only span a small range
along the $r$ direction.  Note also that for such a solution $\int
K\approx\pi$.

\begin{prop}  \label{th:dir}
    (i) A loop with $\int K=2\pi$ reaches the incomplete point from a
    well-defined direction.
    
    (ii) A loop with $\int K>2\pi$ smears out at the incomplete point
    meeting it at an interval of angles.
\end{prop}
\begin{proof}
    Any loop $\gamma_t$ with $\gamma_0$ bounding a disk over which
    $\int K\geq 2\pi$ reaches the incomplete point (in infinite time)
    since it maintains $\int K\geq 2\pi$ throughout the flow so cannot
    disappear.  It is easy to see that (i) implies (ii) since inside a
    loop with $\int K>2\pi$ is a loop with $\int K=2\pi$ symmetric
    around $\theta=\theta_0$ for an interval of $\theta_0$ values, and
    each such loop forms a barrier for the outer loop, forcing it to
    meet the incomplete point at $\theta=\theta_0$.
    
    To prove (i), suppose the converse.  If $\int K=2\pi$ and
    $\gamma_t$ smears out at the incomplete point, take the maximum
    angle $\theta=\theta_{\rm max}$ in the limit of $\gamma_t$. 
    (Recall that all of this takes place in the universal cover, so
    $\theta\in\bbR$.)  This exists, since during the flow, $\gamma_t$
    cannot cross a line $\theta\equiv\theta_1$ so take the minimum of
    the (closed) set of lines $\theta_1$ that lie to the right of
    $\gamma_t$ at some finite time $T$.  Similarly, take the minimum
    angle $\theta_{\rm min}$.
    
    Since $\int K=2\pi$ as $t\rightarrow\infty$, $\gamma_t$ must
    become long and thin in order to reach $\theta_{\rm min}$ and
    $\theta_{\rm max}$ and to maintain its small integral.  Another
    way to state this is to rescale the flow by (\ref{eq:rescale}) to
    get a solution $\tilde{\gamma}_t$ of (\ref{eq:scaleflow}).  The
    static solutions $r=$constant serve as barriers above and below
    the rescaled solution.  The lines $\theta=\theta_{\rm min}$ and
    $\theta=\theta_{\rm max}$ now move outwards exponentially.  The
    curve $\tilde{\gamma}_t$ becomes very flat along some $r=$constant
    in order to maintain $\int K=2\pi$.  It moves outwards
    exponentially fast.  Now, place a static solution as in
    Figure~\ref{fig:scale} with one of its verticals at $(a,b)$ just
    to the right of $\tilde{\gamma}_t$.  The height of
    $\tilde{\gamma}_t$ is much less than that of a static solution
    since by translating the $\theta$ coordinate we can make the
    integral $\int K$ of $\tilde{\gamma}_t$ to the right of the line
    $\theta=0$ arbitrarily small.  The static solution always
    maintains $\int K\approx\pi$ to the right of the line $\theta=0$
    so it is much thicker (spans a greater range along the $r$
    direction) than $\tilde{\gamma}_t$.  Since it is static it will
    serve as a barrier to $\tilde{\gamma}_t$, preventing it from
    moving outwards exponentially as assumed.
    
    The argument does not prevent loops with $\int K>2\pi$ from moving
    outwards exponentially (under the rescaled flow) since the static
    solution is too thin to act as a barrier for taller solutions.
\end{proof}

\subsection{Global flow}

For the remainder of this section we will consider the two-sphere
minus a point equipped with an incomplete metric with the properties:
\begin{equation}  \label{eq:asympos}
    \begin{array}{l}
    \bullet\ {\rm it\ has\ positive\ Gaussian\ curvature}\\
    \bullet\ {\rm it\ is\ asymptotic\ to\ (\ref{eq:modmet})\ or\
	(\ref{eq:cornmet})\ near\ the\ incomplete\ point.}
    \end{array}	
\end{equation}
Over a smooth Riemannian 2-sphere, an embedded closed geodesic
separates it into two pieces with integral of $K$ over each equal to
$2\pi$.  In the incomplete case, we begin with an embedded loop with
integral of $K$ over its interior $2\pi$ and infinite on its exterior
and expect it to flow to an embedded closed geodesic with.

\begin{prop}   \label{th:limits}
    On $S^2-\{\rm point\}$ equipped with a metric satisfying
    (\ref{eq:asympos}), a smooth loop bounding a region with $\int
    K=2\pi$ converges under the flow (\ref{eq:flow}) to one of the
    following:
    
    (i) an embedded closed geodesic;
    
    (ii) the incomplete point, approaching from a well-defined direction;
    
    (iii) the double of a geodesic arc beginning and ending at the
    incomplete point.
\end{prop}
\begin{proof}
    If the curve remains a bounded distance from the incomplete point
    then one can apply the results of \cite{CZhCur}, that the geodesic
    flow produces an embedded closed geodesic, since the proof uses
    only upper and lower bounds for $K$.  Hence (i) is true.
    
    Suppose a loop does not converge to an embedded closed geodesic. 
    Then at least part of the loop must flow arbitrarily close to the
    incomplete point.  If there are disconnected intervals of
    directions at the incomplete point such that the loop moves
    arbitrarily close to the incomplete point from these directions,
    but misses an interval of directions in between, then the number
    of such intervals of directions must be exactly two, the intervals
    must be points, and the loop must be converging to (the double of)
    a geodesic beginning and ending at the incomplete point along the
    two directions.  This is because outside a neighbourhood of the
    incomplete point, arcs of the loop must converge to geodesics.  If
    the limit direction set consists of more than two points, then it
    must contain intervals in order to accommodate different geodesics
    leaving the set.  In that case, the integral $\int K$ over the
    interior of the loop would be infinite.  Alternatively, studying
    scale invariant solutions locally, one can show that the integral
    of the Gaussian curvature on the interior of a solution to the
    flow concentrates near the incomplete point along a well-defined
    direction with $\int K\approx\pi$.  Thus, a loop with $\int K=2\pi$
    can afford only two such directions.  Hence (iii) results.

    If neither (i) nor (iii) results then the loop must be eventually
    contained entirely in an arbitarily small disk neighbourhood of
    the incomplete point.  The neighbourhood of the incomplete point
    is asymptotically given by (\ref{eq:modmet}), (or
    (\ref{eq:cornmet}) near vertices of polygons arising from toric
    surfaces.)  Since the loop moves inside arbitrarily small
    neighbourhoods of the incomplete point, the model metric
    approximates the given metric arbitrarily closely.  By continuity
    of solutions of the flow, we can conclude from the local study
    that the loop takes infinite time to reach the incomplete point,
    that it remains for all time inside any small neighbourhood of the
    incomplete point, and that it approaches it from a well-defined
    direction.
\end{proof}

To complete the proof of Theorem~\ref{th:main} we use the topology of
the disk and its boundary.  First we need the following result about 
small disks.
\begin{lemma}  \label{th:disk}
    Over $S^2-\{\rm point\}$ equipped with a metric satisfying
    (\ref{eq:asympos}) one can continuously assign a family of
    expanding disks around any given point $x\in S^2-\{\rm point\}$ so
    that the integral of $K$ over the disk is arbitrarily large.  The
    assignment is continuous in $x$ and $\int K$.
\end{lemma}
\begin{proof}
    Any well-defined family of disks is suitable for our purposes.  We
    can arbitrarily choose a background metric and take fixed radius
    disks.  In terms of the local picture $ds^2=dr^2+rd\theta^2$ we
    choose a background metric $ds^2=dr^2+d\theta^2$.  We may assume
    the point $x=(a,0)$ by translating $\theta$.  Consider the disk
    $(r-a)^2+\theta^2=a^2$.  Then
    \[\int_DKdA=\int_0^{2a}\int_{-\sqrt{r(2a-r)}}^{\sqrt{r(2a-r)}}
    \frac{r^{1/2}d\theta dr}{4r^2}=\int_0^{2a}\frac{\sqrt{2a-r}}
    {2r}dr=\infty.\] 
    When the metric is only asymptotically like (\ref{eq:modmet}),
    then almost any continuous assignment will do, as long as it
    locally looks like the construction above in a neighbourhood of
    the incomplete point.  In the case of polygons arising from toric
    geometry, it is convenient to choose the disks to be (smoothed)
    polygons centred at each interior point similar to the given
    polygon.
\end{proof}
\begin{cor}  \label{th:diskc}
    One can continuously assign to each point on a two-sphere minus a
    point equipped with a metric satisfying (\ref{eq:asympos}) a disk
    over which $\int K=2\pi$.
\end{cor}
\begin{proof}
    This is immediate from Lemma~\ref{th:disk} since the integral of
    $K$ over a very small disk will be small and a large enough disk
    will give large integral of $K$.
\end{proof}
Before stating the next result, we explain what is meant when ``the
shortest geodesic is long'' which appears in the statement of
Theorem~\ref{th:main}.

The {\em length} of a homology class of loops in $S^2$ is defined as
follows.  Let $\eta$ be a $k$-dimensional family of loops in $S^2$
representing a homology class $[\eta]\in H_k({\rm loops}\subset S^2)$
and take $l(\eta)=$length of the longest loop in $\eta$.  Define
$l([\eta])$ to be the infimum of $l(\eta)$ over all representatives
$\eta$.  If there is an index two (possibly index three when the
condition $\int K=2\pi$ is dropped) embedded closed geodesic on the
incomplete two-sphere, then its length is expected to be the length of
the second homology class of loops.  Note that the assignment of a
family of disks with $\int K=2\pi$ given by Corollary~\ref{th:diskc}
gives an upper bound for the length of the second homology class of
unparametrised embedded loops in $S^2$.  Simply take the maximum
length in the family.  This can be calculated in many cases, in
particular for the metrics arising from the round three-sphere,
$\bbC\bbP^2$ and $S^2\times S^2$.  In these cases, one uses a family
of Euclidean circles, respectively (smoothed) triangles and squares. 
The maximum length loop in the family comes from simple calculus.

The length of the shortest geodesic arc running from the incomplete 
point back to itself can also often be calculated, in particular for 
the same three examples as in the previous paragraph.  In these three 
cases, the geodesic arc can be explicitly calculated.

When twice the length of the shortest geodesic arc running from
the incomplete point back to itself is greater than the length of the
second homology class of embedded loops in the two-sphere we expect 
the geodesic curvature flow to produce an embedded closed geodesic.  
In fact, we need something a little stronger to hold:
\begin{defn}    \label{th:long}
    The shortest geodesic arc running from the incomplete point back
    to itself is {\bf long} if twice its length is greater than the
    maximum length loop in the family constructed in
    Corollary~\ref{th:diskc}.
\end{defn}

\begin{thm}   \label{th:longex} 
    Suppose a metric on a two-sphere minus a point satisfies
    (\ref{eq:asympos}) and its shortest geodesic arc running from the
    incomplete point back to itself is long.  Then there exists a loop
    that converges to an embedded closed geodesic under the geodesic
    curvature flow.
\end{thm}
\begin{proof}
    For each point on the two sphere minus the incomplete point,
    continuously choose the disk centred at the point over which $\int
    K=2\pi$ given by Corollary~\ref{th:diskc}.  The boundary of the
    disk is the initial loop in the geodesic flow.
    
    We will argue by contradiction that one of the loops must flow to
    an embedded closed geodesic.  First, replace the incomplete point
    with its circle of directions, so the two-sphere with an
    incomplete point becomes a disk with circle boundary.  Each loop
    must converge to (i), (ii) or (iii) of
    Proposition~\ref{th:limits}.  No loop can converge to (iii),
    double a geodesic arc, since the length of each initial loop is
    less than twice the length of the shortest geodesic arc, and the
    flow is length decreasing.  Now suppose that no loop flows to (i),
    an embedded closed geodesic.  Thus each loop flows to (ii), the
    incomplete point along a well-defined direction, or equivalently
    to a boundary point on the circle.  This gives a map from the disk
    to its boundary, defined by sending the centre point of an initial
    disk to the destination of its loop on the boundary.
    
    The map is continuous by continuity of the partial differential
    equation governing the flow.  Of course, continuity only applies
    to finite time intervals, but by the proof of
    Proposition~\ref{th:dir}, a loop evolves quite predictably from a
    very large time until infinity, confined to a cone containing the
    well-defined limit direction.  The continuous map can be extended
    to the boundary of the disk by setting it to be the identity. 
    This is because an initial disk centred close to a boundary point
    evolves close to the boundary point by the proof of
    Proposition~\ref{th:dir}.
    
    A continuous map from the closed disk to its boundary that is 
    the identity on the boundary cannot exist by algebraic topology 
    considerations ($S^1\hookrightarrow D^2\rightarrow S^1$ induces 
    the identity on the fundamental group that factors through the 
    zero map) thus giving the desired contradiction.
\end{proof}

The flow simplifies considerably if there is a cyclic symmetry of
order at least three.  In that case, if one begins with a loop that
possesses the symmetry neither (ii) nor (iii) of
Proposition~\ref{th:limits} can occur since if a point of a loop
approaches the incomplete point from some direction then it must
approach from at least three directions, contradicting (ii) and (iii). 
Thus, (i) remains, so the loop must converge to an embedded closed
geodesic.

\section{Minimax and sweepouts.}  \label{sec:minimax}
A sweepout of a manifold is a foliation of the manifold by a path of
codimension 1 submanifolds that degenerate to lower dimensional
submanifolds at each end.  An example of a sweepout is the constant
radius loops on the the two-sphere with a circle invariant metric
described in Section~\ref{sec:circsym}.  The loops degenerate to a
point at each end of the sweepout.  We showed that the maximum length
loop is a geodesic.  For more general sweepouts of a two-sphere by
loops degenerating to points at each end, it is no longer true that
the maximal length loop is a geodesic, however the {\em minimax} in
the family of sweepouts is minimal.  The minimax is obtained from a
family of sweepouts by taking a sequence of sweepouts in the family
with maximal length loop decreasing and converging to the infimum of
all maximal lengths in the family.  A priori the minimax may not be a
smooth loop, and hence a geodesic, however in \cite{PitReg} it is
shown that over a complete Riemannian surface the minimax is a
geodesic.

\subsection{Change of topology in a minimax sequence.}  
\label{sec:change}

Here we give two explicit minimax sequences to show how one can fail
to detect an embedded closed geodesic on an incomplete surface using
sweepouts.  This corresponds in higher dimensions to a change of
topology of the limit from the topology of the leaves of each
sweepout.  The two examples are based on the same idea.

\begin{figure}[ht]
\begin{center}
\scalebox{0.5}{\includegraphics[height=8cm]{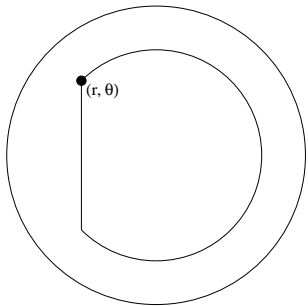}}
\end{center}
\caption{The loop is determined by its left upper corner.}
\label{fig:sweepout}
\end{figure}

Consider the two-sphere with incomplete metric
$ds^2=dr^2+r^2(1-r^2)d\theta^2$.  Any point $(r,\theta)$ with 
$0<\theta<\pi$ determines a loop in the disk given by joining $(r,\theta)$ 
to $(r,-\theta)$ by a vertical line and by a constant radius circle 
to the right of the line as in Figure~\ref{fig:sweepout}.  
The length of the loop is
\[ L=2r\sin\theta\sqrt{1-r^2\cos^2\theta}+2r\theta\sqrt{1-r^2}.\]
(This is clear for $L(\theta=\pi)=2\pi r\sqrt{1-r^2}$,
$L(\theta=\pi/2)=2r+\pi r\sqrt{1-r^2}$, and when $\theta\approx 0$, the
two summands of $L$ are approximately equal.)  A sweepout by loops 
that degenerate to the boundary at one end and to an interior point at 
the other end is simply given by a path $r(\theta)$ such that 
$r(\pi)=1$---the boundary is an end point---$r'(\theta)>0$---the 
leaves are disjoint---and $r(0)\geq 0$---the endpoint is an interior 
point.  

Along $\theta=\pi/2$, the length has a maximum at $r\approx .83$ and
decreases for $r>.83$.  Choose $c>.83$.  Along $r=1$, the length is
increasing as $\theta$ travels from $0$ to $\pi/2$ and decreasing from
$\pi/2$ to $\pi$.  From this it is easy to see that one can take a
path from the point $(c,\pi/2)$ to $(1,\pi)$ such that $L$ decreases
along the path and similarly a path from $(c,\pi/2)$ to $(r_0,0)$ 
along which $L$ increases.  In other words there is a sweepout with 
maximum length loop represented by $(c,\pi/2)$.  As $c\rightarrow 1$ 
the maximum length loop converges to the minimax which is simply the 
line $\theta=\pm\pi/2$.

The above example corresponds to minimal tori in the three-sphere.  We
know that there is an embedded closed geodesic, corresponding to the
Clifford torus, and this is missed by taking the minimax of sweepouts.  
The sweepouts by loops correspond to sweepouts by tori on the 
three-sphere.  The maximal volume tori converge to a two-sphere union 
a one-dimensional arc meeting the two-sphere at two points.  We remove 
the arc, leaving a minimal two-sphere.

The next example uses the metric over the triangle
\begin{equation}   \label{eq:trimet}
    ds^2=y(1-x-y)dx^2+x(1-x-y)dy^2+xyd(1-x-y)^2
\end{equation}
defined in Section~\ref{sec:poly} where we proved the existence of an
embedded closed geodesic.  In the following proposition we show that
this geodesic is missed by taking the minimax of sweepouts.  This
corresponds to sweepouts of tori in $\bbC\bbP^2$ with minimax an
embedded minimal three-sphere.  Inside the maximal volume $T^3$ of
each sweepout, there is a $[0,1]\times T^2\subset T^3$ that collapses
to a pair of disks $[0,1/2]\times S^1\cup[1/2,1]\times
S^1/\{1/2\}\times S^1$.  So the minimax is a three-sphere union two
two-dimensional disks that intersect each other at a point and the
three-sphere at a circle.  We remove the lower-dimensional part of the
limit.

\begin{prop}  \label{th:torisph}
    There exists a sequence of sweepouts of loops of the triangle with
    metric (\ref{eq:trimet}) with minimax an arc that starts and ends
    at the incomplete point.
\end{prop}
\begin{proof}
    Consider sweepouts of the triangle $x\geq 0,\ y\geq 0,\ x+y\leq 1$
    by similar triangles symmetric under $x\leftrightarrow y$.  (There
    are many such sweepouts, each given by a continuous monotone
    function encoded by the top left vertex of each triangle.)
    
    Given any $a>0$ small enough, we will construct a sweepout of the
    triangle by a family of similar triangles such that the maximum
    length triangle is close to the triangle $T_a$ with edges $x=a,\
    y=a,\ x+y=3/4$.  Furthermore, there is a sequence $a_i\rightarrow
    0$, such that the maximum length triangles in the associated
    sweepouts decrease in length.  The minimax of such a sequence is
    the triangle $x=0,\ y=0, x+y=3/4$.  We throw away $x=0$ and $y=0$
    to be left with the geodesic $x+y=3/4$ that meets the boundary.
    
    The length of the path $x+y=c$ running from $x=0$ to $y=0$ is
    $l(c)=c\sqrt{2c(1-c)}$.  By considering $l^2$ it is easy to see
    that the maximum occurs at $c=3/4$ and the function is monotone
    outside the maximum.  The length of the the triangle $T_a$, with
    edges $x=a,\ y=a,\ x+y=3/4$, is
    $L(a)=(3/4-2a)(2\sqrt{2a(1-a)}+\sqrt{3/8})$.  Since
    $L'(a)\rightarrow+\infty$ as $a\rightarrow 0$, for small enough
    $a$, $L(a)>L(0)=l(3/4)>l(c)$ for $c\neq 3/4$.
    
    Given $\epsilon>0$, there is a $\delta>0$ such that
    $l(c)<l(3/4)-\delta$ for $c>3/4+\epsilon$.  Choose $a$ small
    enough so that the path $x=a$ has length less than $\delta/2$. 
    Foliate the outside of the triangle $T_a$ with similar triangles
    (invariant under $x\leftrightarrow y$.)  Then, any triangle with
    side $x+y=c$ and $c>3/4+\epsilon$ has length bounded above by
    $l(c)+2\delta/2<l(3/4)<L(a)$.
    
    Inside $T_a$ we can foliate by triangles each of length less than
    $L(a)$ as follows.  The paths $x+y=c$ for $c<3/4$ running from
    $x=a$ to $y=a$ have lengths less than the path $x+y=3/4$ running
    from $x=a$ to $y=a$ since they are the same proportion of the full
    paths from $x=0$ to $y=0$.  In the extreme, if we were to take the
    triangles with edges $x+y=c$, $x=a$ and $y=a$, for $c$ running
    between $3/4$ and $a$, they would all have length strictly less
    than the triangle $T_a$.  This set of triangles doesn't foliate
    the interior of $T_a$, but we can adjust them slightly to foliate,
    as follows.  Choose the vertical $x=x(c)$ so that the increase in
    the lengths of the path $x=x(c)$ and the path $x=a$ from $y=a$ to
    $y=3/4-x(c)$ (respectively $y=3/4-a$) is exactly equal to the
    decrease in the lengths of $x+y=c$ and $x+y=3/4$ from $x=a$ to
    $y=a$.  (If $x(c)$ gets big enough so that the latter lengths do
    not decrease, then choose any foliation of verticals.)  The
    lengths of all of these triangles are less than the length of
    $T_a$ since we have measured a bit more than necessary.
    
    Thus, the maximum length triangle in the sweepout contains the
    side $x+y=c$ for $c<3/4+\epsilon$.  The maximum length is bounded
    above by $L(a)+\delta$ where $\delta\rightarrow 0$ as
    $\epsilon\rightarrow 0$.  Thus, we can choose a sequence 
    $a_i\rightarrow 0$ so that the maximum length decreases.  
    Moreover, it converges to $L(0)$ and the maximum length triangle 
    converges to $T_0$, the triangle with edges $x=0,\ y=0,\ x+y=3/4$.
\end{proof}

\subsection{Finite dimensions.}
Consider the problem of locating the critical points of a smooth
function $f:M\rightarrow\bbR$ defined on a compact manifold $M$.  The
local minima of $f$ can be found by flowing from a generic point along
the path of steepest descent---the gradient flow, or simply any path
of descent.  A more sophisticated method is needed to locate critical
points of higher index.

In place of a generic point in $M$, take a submanifold $\Sigma\subset
M$, that represents a non-trivial homology class.  If $[\Sigma]\in
H_k(M)$ is non-trivial, then there is a point of $\Sigma$ that flows
to an index $k$ critical point of $M$, under the gradient flow.  This
follows from Morse theory ideas, where the critical points of $f$
represent cohomology classes on $M$ (or more precisely, cochains), and
the non-trivial evaluation of a cohomology class on $[\Sigma]$ detects
the intersection of $\Sigma$ with the stable manifold of a critical
point.

From $\Sigma\subset M$, one might use the gradient flow to locate the 
index $k$ critical point, although this method is limited when we 
adapt it to infinite-dimensional analogues where the gradient flow is 
not so well-behaved.  

Alternatively, a minimax argument can be used.  Take the point in the
image of $\Sigma$ that is a maximum under $f$.  Now vary $\Sigma$ in
its homology class in such a way that its maximum decreases.  The {\em
minimax} of $\Sigma$ is the minimum such maximum value over all
homologous maps $\Sigma\rightarrow M$.  One can prove that it is the
index $k$ critical point.  This is because the maximum point on
$\Sigma$ becomes arbitarily close to the intersection of $\Sigma$ with 
the stable manifold of the critical point, and this intersection can 
be brought arbitrarily close to the critical point.

\begin{figure}[ht]
\begin{center}
\scalebox{0.5}{\includegraphics[height=8cm]{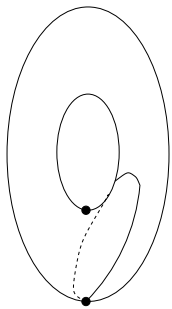}}
\end{center}
\caption{The index 1 critical point is a minimax.}
\label{fig:minimax}
\end{figure}

An example is given in Figure~\ref{fig:minimax} for the height
function on a torus.  Consider maps $\gamma:[0,1]\rightarrow M$ that
satisfy $\gamma(0)=q=\gamma(1)$ where $q\in M$ is the minimum, and the
image of $\gamma$ represents a non-trivial class in $H_1(M)$.  The
example shows that the image of $\gamma$ contains a point that flows
to the index 1 critical point under the gradient flow.  The maximum
point on $\gamma$ is forced to lie above the index 1 critical point,
and a minimax sequence of such maximum points will converge to the
index 1 critical point.

This shows that sweepouts can give index one closed geodesics at best
and that we need to take the minimax of higher dimensional families of
loops to get some closed geodesics.  Lusternik-Schnirelmann theory
successfully does this on the infinite-dimenional manifold of loops on
a manifold equipped with the length functional.  More generally,
Almgren \cite{AlmGeo} used geometric measure theory techniques to
prove the existence of minimax limits and Pitts \cite{PitReg} proved
regularity results for the varifolds produced from Almgren's work.

One way of producing the minimax of a family of paths on manifold is
to apply a curve shortening flow to the entire family.  As long as
such a process is continuous, the topology of the underlying manifold
can be used to deduce the existence of closed geodesics.  Hass and
Scott \cite{HScSho} constructed a curve shortening flow by covering a
surface with small disks and then straightening a given path in each
disk, in turn.  They proved that this process produces a sequence of
paths that converge to a geodesic.  Furthermore, any given continuous
family of paths can be shortened continuously---this requires some
fixing of discontinuities.  When the surface has incomplete points,
this process can still be used to shorten a given curve.  However, the
curve shortening may not be unique and this produces discontinuities
in the curve shortening of a family that cannot be fixed. 
\begin{figure}[ht]
\begin{center}
\scalebox{0.8}{\includegraphics[height=8cm]{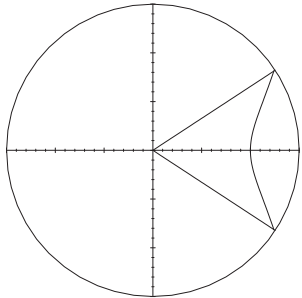}}
\end{center}
\caption{Equal length shortest paths joining boundary points.}
\label{fig:nonunique}
\end{figure}
Figure~\ref{fig:nonunique} shows two shortest paths of the same length
joining a pair of points on the boundary of a small disk around the
incomplete point with metric $ds^2=dr^2+rd\theta^2$.  Arbitrarily
small disks possess such pairs since in a disk of radius $\epsilon$
geodesics (with $\theta$ non-constant) from the boundary that come
close enough to the incomplete point have length greater than
$2\epsilon$ so the intermediate value theorem produces a geodesic from
the boundary with length exactly $2\epsilon$ and it can be shown it is
minimal.  This results in a choice for the curve shortening flow,
which produces a discontinuity in the flow of a family.  A similar
type of discontinuity occurs in \cite{HScSho}, and there they fill in
the ``gap'' between the two shortest paths with a family of paths.  In
our case, the family of paths joining the two boundary points consists
of paths longer than the two bounding shortest paths, and this
destroys the {\em shortening} of a family.  The key difference with
\cite{HScSho} is that their discontinuity arises between a shortest
path and a (local) longest path.

\section{Minimal submanifolds.}
A natural source of incomplete metrics over a surface comes from the
study of minimal submanifolds in higher dimensions.  If $Y\rightarrow
M$ is a codimension one minimal immersion invariant under a subgroup
$G$ of isometries of $M$ such that $\dim M/G=2$, then it is the
pull-back of a geodesic on the quotient with an adjusted metric.  If
$g$ is the quotient metric on $M/G$, and $V(x)$ is the volume of the
orbit of $x\in M$ under $G$ (defined with multiplicity), then the
adjusted metric is $V^2g$ \cite{HLaMin}.  Since the volume of each 
orbit is taken into account, the length of a geodesic in the quotient 
is equal to the volume of its pull-back upstairs.

\subsection{Spheres}

Consider $S^3$ with the circle action that fixes a circle.  The quotient
space of the action is a disk with boundary corresponding to the fixed
point set upstairs.  The lengths of the orbits go to zero as the fixed
point set is approached so this produces a metric on the disk which is
degenerate on the boundary, or equivalently an incomplete metric on
the two-sphere minus a point.  For concreteness, we will consider the
special case of circle invariant metrics on $S^3$ of the form 
\[ ds^2=dr^2+u'(r)^2d\theta^2+\frac{f(u(r),\phi)^2}{u'(r)^2}d\phi^2\]
where $f(u)^2\sim u^2(1-u^2)$ near $u=0$ and $u=1$ and
$u(r)\sim\sin(r)$ near $r=0$ and $r=\pi/2$ (strictly we mean
$u'(0)=1$, etc.)  This means the metric is asymptotically the same as
the round metric at the two invariant circles $\theta=0$ and $\phi=0$. 
The more general metric would still look like this asymptotically,
although it would not be diagonal, and the coefficient of $d\theta^2$
would depend on $\phi$ too.  The quotient metric is
\[ d\bar{s}^2=dr^2+\frac{f(u(r),\phi)^2}{u'(r)^2}d\phi^2\]
and the adjusted metric is
\[ d\hat{s}^2=u'(r)^2dr^2+f(u(r),\phi)^2d\phi^2=du^2+f(u,\phi)^2d\phi^2.\]
This is complete at $u=0$ and incomplete at $u=1$ with $\int K=\infty$ 
in a neighbourhood.

More generally, the $SO(n-1)$ action on the first factor of
$\bbR^{n-1}\times\bbR$ induces an action on $S^n$.  The quotient is a
disk with boundary corresponding to the fixed point set upstairs, so
again we get an incomplete metric on the two-sphere minus a point. 
Near the incomplete point, the metric looks like
\[ ds^2=dr^2+r^pd\theta^2,\ p=2-2/n.\]

When 
\begin{equation}  \label{eq:round}
    ds^2=du^2+u^2(1-u^2)d\theta^2 
\end{equation}    
corresponding to the round three-sphere, the length of the second
homology class of loops with $\int K=2\pi$ over the disks they bound,
is the same as the length of the embedded closed geodesic which is
$\pi$.  The length of the shortest geodesic arc is 2.  Since
$\pi<2.2$, the length of the second homology class of loops is less
than twice the length of the shortest geodesic arc.  The same is true
for a metric sufficiently close to (\ref{eq:round}), and furthermore a
sufficiently close metric has positive Gaussian curvature.  Thus we
have proven Theorem~\ref{th:sphere}.  Note that, as mentioned in the
introduction, a three-sphere that gives rise to a positive Gaussian
curvature incomplete two-sphere may be far away from the round
three-sphere if it also possesses a cyclic symmetry of order at least
three.

\subsection{Toric varieties}   \label{sec:toric}

A Kahler four-manifold is {\em toric} if it admits an action of the
torus $T^2$ that preserves the Kahler structure.  We will denote it by
its underlying symplectic manifold $(M^4,\omega)$ and prove results
for a large class of compatible toric Kahler metrics.  It is fibred by
Lagrangian tori over a convex polygon base $P\subset\bbR^2$ with some
degenerate fibres.
\[\begin{array}{ccc}T^2&\hookrightarrow&(M^4,\omega)\\
&&\downarrow\\&&\ \ \ \ \ \ \ \ P\subset\bbR^2\end{array}.\] 
The quotient equipped with the adjusted metric of Hsiang and Lawson 
is an incomplete metric on the two-sphere minus a point.  Before we 
study this specific case of four-dimensional toric manifolds we will 
describe some background in any dimension.

Let $(M^{2n},\omega)$ be a closed symplectic manifold equipped with an
effective action of a torus $T^n$ that preserves $\omega$.  A basic
example of this is complex projective space, $\bbC\bbP^n$, with its
natural symplectic form $\partial\bar{\partial}\ln|z|^2$ where
$z=(z_0,..,z_n)$ is a projective coordinate.  It is invariant under
the action of $U(n+1)$, and in particular under $T^{n+1}\subset
U(n+1)$ which acts as $T^n$ (since one circle acts trivially.) 
Another basic example is the subset of all Hermitian matrices with a
given set of eigenvalues.  The torus acts by conjugation and the
symplectic form uses the trace and commutator of matrices.  (The first
example is a case of the latter: the space of $(n+1)\times(n+1)$
Hermitian matrices with one eigenvalue $1$ and $n$ eigenvalues $0$, is
isomorphic, as a symplectic manifold with torus action, to
$\bbC\bbP^n$.)

Such an action is said to be {\em Hamiltonian}, and every $\xi\in\bf
t=\bbR^n$ gives rise to the function $H_{\xi}:M\rightarrow\bbR$ that
satisfies $dH_{\xi}=\omega(\tilde{\xi},\cdot)$ for $\tilde{\xi}$ the
vector field generated by $\xi$ and the torus action.  Define the {\em
moment map} to be $\phi:M\rightarrow\bf t^*=\bbR^n$ that satisfies
\[(\phi(p),\xi)=H_{\xi}(p).\] In the case of $\bbC\bbP^n$,
\[\phi(z_0,..,z_n)=\left(\frac{|z_1|^2}{|z|^2},..,\frac{|z_n|^2}{|z|^2}\right)\]
which has image $\{(x_1,..,x_n)|x_i\geq 0,\ \sum x_i\leq 1\}$.  In
general, the image of $\phi$ is a convex polytope $P$ in $\bbR^n$. 
This generalises the theorem of Schur that the diagonal elements of an
$n\times n$ Hermitian matrix lie in the convex hull of the points in
$\bbR^n$ obtained from the eigenvalues of the matrix, listed in every
possible order.

Those polytopes that arise as the image of a moment map of a torus
action are characterised in the following definition.
\begin{defn}
    A convex polytope $P\subset\bbR^n$ is Delzant if:
    
    (1) there are $n$ edges meeting at each vertex $v$;
    
    (2) the edges meeting at the vertex $v$ are rational, i.e. 
    $v+tv_i,0\leq t\leq\infty$, $v_i\in\bbZ^n$;
    
    (3) the $v_1,\ldots,v_n$ can be chosen to be a generating set of 
    $\bbZ^n$.
\end{defn}
The Delzant polytope $P$ can be defined by inequalities
\begin{equation}  \label{eq:poly}
    l_r(x)=\langle x,\mu_r\rangle-\lambda_r\geq 0
\end{equation}    
where $\mu_r$ is a primitive element of $\bbZ^n$ giving the inward
pointing normal of the $r$th face of $P$.  For each Delzant polytope
$P$, Guillemin constructed a symplectic manifold $M_P$ using the
beautiful idea of varying the numbers $\lambda_r$.  The variations
generate a torus action, and one obtains $M_P$ as a symplectic
quotient of this action.  See \cite{GuiKah} for details.  There is a
bijective correspondence between Delzant polytopes and toric manifolds
$P\mapsto M_P$.
\begin{other}[\cite{DelHam}]
    Let $(M^{2n},\omega)$ be a closed symplectic manifold equipped
    with an effective action of a torus $T^n$ that preserves $\omega$
    and moment map $\phi:M\rightarrow\bf t^*=\bbR^n$.  Then the image
    of $\phi$ is a Delzant polytope, and $M$ and the $T^n$ action are
    isomorphic to $M_P$.
\end{other}

Toric manifolds are Kahler manifolds and so far we have only described
a symplectic structure.  When we also specify a (compatible) complex
structure, this together with the symplectic structure gives rise to
the Kahler metric.

The metric is usually expressed with respect to one of two {\em
canonical coordinate systems}.  Denote by $M^{\circ}\subset M$ those
points on which the torus action is free.  The torus action
complexifies and allows us to describe $M^{\circ}$ in complex
coordinates by 
\[ M^{\circ}\cong\bbC^n/2\pi i\bbZ^n=\bbR^n\times
iT^n=\{u+iv:u\in\bbR^n,v\in\bbR^n/\bbZ^n\}.\] 
Alternatively, one can use action-angle coordinates: 
\begin{equation}  \label{eq:actang}
    M^{\circ}\cong P^{\circ}\times T^n=\{(x,y):x\in
    P^{\circ}\subset\bbR^n,y\in\bbR^n/\bbZ^n\}
\end{equation}
where $P^{\circ}$ is the interior of the Delzant polytope.    

We are interested in the quotient manifold, which is identified with
either $\bbR^n$ or the Delzant polytope $P$, equipped with the
quotient metric.  In fact, the metric on the toric manifold is
expressed in terms of functions over the quotient, $\bbR^n$ or $P$,
explicitly using the quotient metric.  

For our purposes, the action-angle coordinate system is more suitable,
with its compact quotient.  When we adjust the quotient metric by the
volumes of the fibres, the new metric is degenerate at the faces of
the Delzant polytope.

We will follow Abreu's treatment \cite{AbrKah} of Guillemin's
construction \cite{GuiKah} of metrics over toric manifolds.  With
respect to the coordinates (\ref{eq:actang}), the symplectic form is
the standard $\omega=\sum_j dx_j\wedge dy_j$, given in matrix form by
\[\left(\begin{array}{ccc}
0&\vdots&I\\\ldots&\ldots&\ldots\\-I&\vdots&0 \end{array}\right).\]
The complex Kahler structure $J$ is given by
\[\left(\begin{array}{ccc}
0&\vdots&-G^{-1}\\\ldots&\ldots&\ldots\\G&\vdots&0
\end{array}\right)\] 
where $G={\rm Hess}(g)$ is the Hessian of a potential $g(x)$,
\[G=\left[\frac{\partial^2g}{\partial x_j\partial x_k}\right],\ 1\leq
j,k\leq n.\] 

We can now describe the potential $g(x)$ defined over the polytope 
$P$, using the defining expressions for $P$ (\ref{eq:poly}).  Define
\[ g_P(x)=\frac{1}{2}\sum^d_{r=1}l_r(x)\log l_r(x).\]
This is well-defined on the interior of $P$ since $l_r(x)>0$ there.
More generally, choose any function $h$ smooth on the whole of $P$, 
and set $g=g_P+h$.  It follows that 
\begin{equation}  \label{eq:metric}
    \det(G)=\left(\delta(x)\prod^d_{r=1}l_r(x)\right)^{-1}
\end{equation}
for $\delta(x)$ a strictly positive smooth function on $P$.  

The quotient metric on $P$ is given by the $n\times n$ matrix $G$.  It
is a well-defined finite metric at the boundary of $P$ although it
looks like it blows-up there.  For example, on $\bbC\bbP^2$ the
Fubini-Study metric is described by the canonical potential on the
Delzant polygon $x_1\geq 0$, $x_2\geq 0$ and $x_1+x_2\leq 1$,
$g_P=x_1\ln x_1+x_2\ln x_2+(1-x_1-x_2)\ln(1-x_1-x_2)$.  The quotient
metric is the round metric on an octant of the two-sphere (identified
with a triangle) which is certainly finite at the boundary of the
octant.

We can associate minimal surfaces in the quotient to minimal surfaces
in the toric manifold upstairs once we adjust the quotient metric by
the volumes of the torus fibres.  From (\ref{eq:metric}), the volume
$V(x)$ of the torus fibre over $x\in P$ is given by $V(x)^2=\det
G^{-1}(x)$.  

Replace the quotient metric $G$ by 
\[ g(x)=V(x)^{2/(n-1)}G(x)=\det G^{1/(1-n)}(x)G(x).\]
By (\ref{eq:metric}), this satisfies $\det 
g(x)=\left(\delta(x)\prod^d_{r=1}l_r(x)\right)^{1/(n-1)}$
which vanishes on the boundary of $P$.  When $n=2$, it vanishes 
linearly and explicitly when $l_r=a_rx+b_ry+c_r$ we get:
\[ g(x)=\sum\frac{1}{2l_r}\left(\begin{array}{cc}a_r^2&
a_rb_r\\ a_rb_r&b_r^2\end{array}\right)\]
which has determinant $\delta(x)/(\prod l_r)$ where $\delta(x)\neq 0$ 
on the closed polygon.  Equivalently, the quotient metric is
\[ d\bar{s}^2=\frac{1}{2}\sum\frac{dl_r^2}{l_r}\] 
so the adjusted metric is
\[ ds^2=\frac{1}{2\delta}\sum\prod_{j\neq r}l_jdl_r^2.\]
At each edge $l_r=0$ the metric is asymptotic to (\ref{eq:modmet}).

Large families of metrics can be obtained for a given polygon this
way.  Simply add ``edges'' to the polygon that lie outside the
polygon.  For example, over $\bbC\bbP^2$ with polygon defined by
$x\geq 0$, $y\geq 0$ and $1-x-y\geq 0$ we can include the inequality
$l(x,y)=2-x-y\geq 0$ and still get a positive definite metric on the
interior of the triangle.  We may also take multiples $\epsilon
l(x,y)$ for an extra edge, which is desirable when we choose $\epsilon$
small so that the new metric is close to the original metric.  (In
fact one can take $g_P(x)+h(x)$ for any smooth function $h(x)$ defined
on a neighbourhood of $P$, and when $h(x)$ is small enough the new 
metric remains positive definite.)

For the Fubini-Study metric on $\bbC\bbP^2$ we proved the existence of
an embedded closed geodesic.  An easier calculation shows that the
shortest arc is given by the line $x=1/4$ which has length
$3\sqrt{3}/8\sqrt{2}\approx .46$.  The maximum length triangle with
$\int K=2\pi$ on its interior is less than .9 which is less than
$2\times 3\sqrt{3}/8\sqrt{2}$.  Thus, the shortest geodesic is long and we
can apply Theorem~\ref{th:longex} to get Theorem~\ref{th:cp2}.  The
same reasoning applies to $S^2\times S^2$ with the product of round
metrics of the same area.

The Fubini-Study metric on $\bbC\bbP^2$ is invariant under a $\bbZ^3$
action.  We can drop the condition that a metric on $\bbC\bbP^2$ be
sufficiently close to the Fubini-Study metric and instead require that
it is invariant under the $\bbZ^3$ action and that the adjusted metric
on the triangle has positive Gaussian curvature to again deduce the 
existence of an embedded closed geodesic.

The geodesic curvature flow proof corresponds to mean curvature flow
in $\bbC\bbP^2$.  The {\em mean curvature} of a hypersurface
$Y^3\subset M^4$ is a function on $Y^3$ given by the trace of the
second fundamental form of $Y^3$.  As a (symmetric) bilinear form, the
second fundamental form is defined by $l(X,Y)=\langle
\nabla_X\nu,Y\rangle$ for $X,Y\in T_pY^3$, $\nu$ the unit normal
vector field to $Y^3$ and $\nabla$ the Levi-Civita connection on
$M^4$.

The mean curvature flow of a surface $Y$ is an evolution of $Y$ in its
normal direction with magnitude given by the mean curvature: 
\[\frac{dY_t}{dt}={\rm mean\ curvature}\cdot\nu.\] 
In actual fact, we need to adjust the flow further, by multiplying the
right-hand side by a function defined on $M^4$ (and hence independent
of the way in which $Y$ embeds.)  For a hypersurface invariant under 
the torus action we have
\[ {\rm mean\ curvature\ =\frac{1}{\rm area\ of\ torus}}k\]
where $k$ is the geodesic curvature of the corresponding curve on the
polygon.  The reciprocal of the Gaussian curvature of the polygon is
an invariant function $\psi$ on the polygon that seems not to have a
natural interpretation.  Thus the flow in the toric surface is given 
by
\[\frac{dY_t}{dt}=\psi\cdot({\rm area\ of\ torus})^2\cdot
{\rm mean\ curvature}\cdot\nu\]
where an extra factor of the area of the torus arises from comparing
the normal vectors in two-dimensions and four-dimensions.  The factor
$\psi\cdot({\rm area\ of\ torus})^2$ is a function on the toric
surface that vanishes on the divisor with isotropy group.

\section{Counterexample}  \label{sec:countex}
It is not always true that there exists an embedded closed geodesic 
loop on an incomplete two-sphere.  A counterexample is as follows.

Take the metric 
\begin{equation}  \label{eq:neck}
    ds^2=\sin^2(y)(dx^2+dy^2),\ (x,y)\in[-R,R]\times (0,\pi). 
\end{equation}
Choose R very large and cap off the metric at the two ends in such a
way that any geodesic in the cap must leave the cap.  (It can be
chosen to have positive Gaussian curvature.)  For example, choose
\begin{equation}  \label{eq:cap}
    ds^2=\cos^2(|x|-R)\sin^2y(dx^2+\cos^2(|x|-R)dy^2),\ |x|\in[R,R+\pi/2).
\end{equation}
Alternatively, the metric can be capped off smoothly.

\begin{thm}
    There is no embedded closed geodesic loop for the metric 
    (\ref{eq:neck}), (\ref{eq:cap}).
\end{thm}
\begin{proof}
    We call a geodesic vertical if it is given by $x=$constant for
    $x\in[-R,R]$.  The theorem follows from the following facts:
    
    (i) any non-vertical geodesic cannot be vertical anywhere in 
    $x\in[-R,R]$;
    
    (ii) any two non-vertical geodesics must intersect;
    
    (iii) a geodesic loop must enter $x\in[-R,R]$ twice and hence it
    must be self-intersecting.\\
    
    Uniqueness of geodesics gives (i) immediately.  For (ii) we use
    the fact that the geodesic flow is integrable and translation
    invariant in $x\in[-R,R]$, as dealt with in
    Section~\ref{sec:circsym}, to deduce that inside $x\in[-R,R]$ each
    geodesic oscillates around $y=\pi/2$ with period
    \begin{eqnarray*}
	\Omega_c&=&4\int^{\pi/2}_{\sin x=c}\frac{c}{\sqrt{\sin^2x-c^2}}
	dx\\
	&=&4c\int^1_c\frac{du}{\sqrt{(u^2-c^2)(1-u^2)}}\\
	&<&4\sqrt{\frac{c}{2(1+c)}}\int^1_c\frac{du}{\sqrt{(u-c)(1-u)}}\\
	&=&4\pi\sqrt{\frac{c}{2(1+c)}}.
    \end{eqnarray*}
    where $c$ denotes the geodesic that reaches a maximum value of
    $y=c^2$ in $x\in[-R,R]$.  The period is bounded and we choose $R$
    to be greater than (half) this period so that each non-vertical
    geodesic must meet the line $y=\pi/2$.  Since any geodesic meets
    $y=\pi/2$ with a given period, any two non-vertical geodesics meet
    in $x\in[-R,R]$.
    
    Now a geodesic in the cap (that does not meet the boundary) must
    leave the cap since the equation for geodesics in the cap is
    \begin{eqnarray*}
	\ddot{x}&=&\dot{x}^2\tan(x-R)-2\dot{x}\dot{y}\cot y-\dot{y}^2\sin 
	2(x-R)\\
	\ddot{y}&=&\dot{x}^2\sec^2(x-R)\cot y+4\dot{x}\dot{y}\tan(x-R)
	-\dot{y}^2\cot y
    \end{eqnarray*}
    Thus $\dot{x}=0\ \Rightarrow\ \ddot{x}<0$ so no loops can occur 
    inside the cap.  By (ii), the two arcs of any closed geodesic loop 
    must intersect in $x\in[-R,R]$ so (iii) follows and the theorem is 
    proven.
\end{proof}

Geodesics on the disk with metric (\ref{eq:neck}), (\ref{eq:cap})
arise from minimal surfaces on $S^3$ equipped with the circle
invariant metric obtained by stretching the round metric in a
neighbourhood of a minimal two-sphere, or more precisely
$S^2\times\bbR$ with the round metric times the flat metric capped off
with round three balls.  The counterexample shows that there is no
minimally embedded $T^2\subset S^3$ invariant under the circle 
action.  It would be interesting to know if there is any minimally 
embedded $T^2$ in this three-sphere.

\begin{cor}
    There is a length decreasing sequence of loops that converge to a
    double geodesic arc.
\end{cor}
\begin{proof}
    Take a loop that bounds a region with $\int K=2\pi$ (where
    $K=$Gaussian curvature) that is contained inside $x\in[-R,R]$ and
    that is symmetric under reflection in $y=\pi/2$.  Under the flow
    (\ref{eq:flow}) the loop cannot move past vertical geodesics and
    $\int K=2\pi$ is preserved.  It cannot converge to a geodesic loop
    since no vertical tangencies are allowed so it travels to both
    boundaries symmetrically, converging to a double vertical
    geodesic.  The flow is length decreasing so the corollary follows.
\end{proof}

\end{document}